\documentclass[10pt]{article}

\usepackage{amscd,amsmath, amssymb, fancyhdr, mathbbol}

\numberwithin{equation}{section}

\newcommand{\version}{version 6.5,\ \   Feb. 4, 2014}

\def\eqref#1{(\ref{#1})}

\newcommand{\arrow}{{\:\longrightarrow\:}}
\newcommand{\Z}{{\Bbb Z}}
\def\C{{\Bbb C}}
\newcommand{\R}{{\Bbb R}}
\newcommand{\Q}{{\Bbb Q}}

\def\1{\sqrt{-1}\:}
\newcommand{\restrict}[1]{{\left|_{{\phantom{|}\!\!}_{#1}}\right.}}
\newcommand{\cntrct}                % contraction with a vector field
{\hspace{2pt}\raisebox{1pt}{\text{$\lrcorner$}}\hspace{2pt}}

\newcommand{\calo}{{\cal O}}

\renewcommand{\tilde}{\widetilde}
\renewcommand{\bar}{\overline}
\renewcommand{\phi}{\varphi}
\renewcommand{\epsilon}{\varepsilon}
\renewcommand{\geq}{\geqslant}
\renewcommand{\leq}{\leqslant}
\renewcommand{\max}{{\rm max}}

\newcommand{\Pic}{\operatorname{Pic}}

\newcommand{\dist}{\operatorname{ \text {\it \sffamily dist}}}

\newcounter{Mycounter}[section]
\newcounter{lemma}[section]
\newcounter{claim}[section]
\renewcommand{\theclaim}{{Claim \thesection.\arabic{claim}}}
\newcommand{\claim}{%
    \setcounter{claim}{\value{Mycounter}}
    \refstepcounter{claim}
    \stepcounter{Mycounter}
    {\noindent \bf \theclaim:\ }}

\newcounter{sublemma}[section]
\newcounter{corollary}[section]
\renewcommand{\thecorollary}{{Corollary \thesection.\arabic{corollary}}}
\newcommand{\corollary}{%
    \setcounter{corollary}{\value{Mycounter}}
    \refstepcounter{corollary}
    \stepcounter{Mycounter}
    {\noindent \bf \thecorollary:\ }}

\newcounter{theorem}[section]
\renewcommand{\thetheorem}{{Theorem \thesection.\arabic{theorem}}}
\newcommand{\theorem}{%
    \setcounter{theorem}{\value{Mycounter}}
    \refstepcounter{theorem}
    \stepcounter{Mycounter}
    {\noindent \bf \thetheorem:\ }}

\newcounter{conjecture}[section]
\renewcommand{\theconjecture}{{Conjecture \thesection.\arabic{conjecture}}}
\newcommand{\conjecture}{%
    \setcounter{conjecture}{\value{Mycounter}}
    \refstepcounter{conjecture}
    \stepcounter{Mycounter}
    {\noindent \bf \theconjecture:\ }}

\newcounter{proposition}[section]
\renewcommand{\theproposition}
      {{Proposition \thesection.\arabic{proposition}}}
\newcommand{\proposition}{%
    \setcounter{proposition}{\value{Mycounter}}
    \refstepcounter{proposition}
    \stepcounter{Mycounter}
    {\noindent \bf \theproposition:\ }}

\newcounter{definition}[section]
\renewcommand{\thedefinition}
      {{Definition~\thesection.\arabic{definition}}}
\newcommand{\definition}{%
    \setcounter{definition}{\value{Mycounter}}
    \refstepcounter{definition}
    \stepcounter{Mycounter}
    {\noindent \bf \thedefinition:\ }}

\newcounter{example}[section]
\renewcommand{\theexample}{{Example \thesection.\arabic{example}}}
\newcommand{\example}{%
    \setcounter{example}{\value{Mycounter}}
    \refstepcounter{example}
    \stepcounter{Mycounter}
    {\noindent \bf \theexample:\ }}

\newcounter{remark}[section]
\renewcommand{\theremark}{{Remark \thesection.\arabic{remark}}}
\newcommand{\remark}{%
    \setcounter{remark}{\value{Mycounter}}
    \refstepcounter{remark}
    \stepcounter{Mycounter}
    {\noindent \bf \theremark:\ }}

\newcounter{problem}[section]
\newcounter{question}[section]
\newcommand{\proof}{\noindent{\bf Proof:\ }}

\makeatletter

\@addtoreset{equation}{section} \@addtoreset{footnote}{section}
\makeatother

\def\blacksquare{\hbox{\vrule width 5pt height 5pt depth 0pt}}
\def\endproof{\blacksquare}

\begin{document}
%%%%%%%%%%%%%%%%%%%%%%%%%%%%%%%%%%%%%%%%%%%%%%%%%%%%%%%%%%%%
\begin{center}
{\LARGE\bf
Parabolic nef currents\\[2mm] on hyperk\"ahler manifolds\\[4mm]
}
%%%%%%%%%%%%%%%%%%%%%%%%%%%%%%%%%%%%%%%%%%%%%%%%%%%%%%%%%%%%

 Misha
Verbitsky\footnote{Partially supported by RFBR grants
 12-01-00944-Á,  10-01-93113-NCNIL-a, and
AG Laboratory NRI-HSE, RF government grant, ag. 11.G34.31.0023.}

\end{center}

%%%%%%%%%%%%%%%%%%%%%%%%%%%%%%%%%%%%%%%%%%%%%%%%
{\small \hspace{0.1\linewidth}
\begin{minipage}[t]{0.8\linewidth}
{\bf Abstract} \\
Let $M$ be a compact, holomorphically symplectic K\"ahler
manifold, and $\eta$ a (1,1)-current which is nef
(a limit of K\"ahler forms). Assume that the
cohomology class of $\eta$ is parabolic, that is,
its top power vanishes. We prove that all 
Lelong sets of $\eta$ are coisotropic.
When $M$ is generic, this is used to show 
that all Lelong numbers of $\eta$ vanish.
We prove that any hyperk\"ahler
manifold with $\Pic(M)=\Z$ has 
non-trivial coisotropic subvarieties,
if a generator of $\Pic(M)$ is parabolic.
\end{minipage}
}
%%%%%%%%%%%%%%%%%%%%%%%%%%%%%%%%%%%%%%%%%%%%%%%%

{\scriptsize
\tableofcontents
}

%%%%%%%%%%%%%%%%%%%%%%%%%%%%%%%%%%%%%%%%%%%%%%%%

\section{Introduction}

%%%%%%%%%%%%%%%%%%%%%%%%%%%%%%%%%%%%%%%%%%%%%%%%

%%%%%%%%%%%%%%%%%%%%%%%%%%%%%%%%%%%%%%%%%%%%%%%%
\subsection{Hyperk\"ahler manifolds}
%%%%%%%%%%%%%%%%%%%%%%%%%%%%%%%%%%%%%%%%%%%%%%%%

\definition  A {\bf hyperk\"ahler manifold}
is a compact, K\"ahler, holomorphically symplectic manifold.

\hfill

\definition
 A hyperk\"ahler manifold $M$ is called
{\bf simple} if $H^1(M)=0$, $H^{2,0}(M)=\C$.

\hfill

\theorem
(Bogomolov's Decomposition Theorem,
\cite{_Bogomolov:decompo_}, \cite{_Besse:Einst_Manifo_}). 
Any hyperk\"ahler manifold admits a finite covering
which is a product of a torus and several 
simple hyperk\"ahler manifolds. 
\endproof

\hfill

\remark
Further on, all hyperk\"ahler manifolds
are silently assumed to be simple.

\hfill

{\bf A note on terminology.}
Speaking of hyperk\"ahler manifolds, people
usually mean one of two different notions.
One either speaks of holomorphically
symplectic K\"ahler manifold, or of 
a manifold with a {\em hyperk\"ahler structure},
that is, a triple of complex structures
satisfying quaternionic relations and
parallel with respect to the Levi-Civita
connection. The equivalence
(in compact case) between these two 
notions is provided by the Yau's solution
of Calabi-Yau conjecture 
(\cite{_Besse:Einst_Manifo_}). Throughout this paper,
we use the complex algebraic geometry
point of view, where ``hyperk\"ahler''
is synonymous with ``K\"ahler 
holomorphically symplectic'', 
in lieu of the differential-geometric
approach. To avoid the terminological
confusion, we tried not mention quaternionic 
structures (except Subsection 
\ref{_gene_hyper_subva_Subsection_}, where it was 
impossible to avoid). 
The reader may
check \cite{_Besse:Einst_Manifo_} for an introduction
to hyperk\"ahler geometry from the
differential-geometric point of view.

Notice also that we included compactness in our definition
of a hyperk\"ahler manifold. In the differential-geometric 
setting, one does not usually assume that the
manifold is compact.

%%%%%%%%%%%%%%%%%%%%%%%%%%%%%%%%%%%%%%%%%%%%%%%%%%%%%%%%%%%%
\subsection{The Bogomolov-Beauville-Fujiki form}
%%%%%%%%%%%%%%%%%%%%%%%%%%%%%%%%%%%%%%%%%%%%%%%%%%%%%%%%%%%%

\theorem
(\cite{_Fujiki:HK_})
Let $\eta\in H^2(M)$, and $\dim M=2n$, where $M$ is
hyperk\"ahler. Then $\int_M \eta^{2n}=\lambda q(\eta,\eta)^n$,
for some primitive integer quadratic form $q$ on $H^2(M)$
and $\lambda \in \Q^{>0}$.
\endproof

\hfill

\definition
This form is called
{\bf  Bogomolov-Beauville-Fujiki form}.  It is defined
by this relation uniquely, up to a sign. The sign is determined
from the following formula (Bogomolov, Beauville;
\cite{_Beauville_}, \cite{_Huybrechts:lec_}, 23.5)
\begin{align*}   \lambda q(\eta,\eta) &=
   (n/2)\int_X \eta\wedge\eta  \wedge \Omega^{n-1}
   \wedge \bar \Omega^{n-1} -\\
 &-(1-n)\frac{\left(\int_X \eta \wedge \Omega^{n-1}\wedge \bar
   \Omega^{n}\right) \left(\int_X \eta \wedge
   \Omega^{n}\wedge \bar \Omega^{n-1}\right)}{\int_M \Omega^{n}
   \wedge \bar \Omega^{n}}
\end{align*}
where $\Omega$ is the holomorphic symplectic form,
and $\lambda$ a positive constant.

\hfill

\remark
The form $q$ has signature $(b_2-3,3)$.
It is negative definite on primitive forms, and positive
definite on the space $\langle \Omega, \bar \Omega, \omega\rangle$
 where $\omega$ is a K\"ahler form, as seen from the
following formula
\begin{multline}\label{_BBF_via_Kahler_Equation_}
   \mu q(\eta_1,\eta_2)= \\
   \int_X \omega^{2n-2}\wedge \eta_1\wedge\eta_2  
   - \frac{2n-2}{(2n-1)^2}
   \frac{\int_X \omega^{2n-1}\wedge\eta_1 \cdot
   \int_X\omega^{2n-1}\wedge\eta_2}{\int_M\omega^{2n}}, \
   \  \mu>0
\end{multline}
(see e. g. \cite{_Verbitsky:cohomo_}, Theorem 6.1,
or \cite{_Huybrechts:lec_}, Corollary 23.9).

\hfill

%%%%%%%%%%%%%%%%%%%%%%%%%%%%%%%%%%%%%%%%%%%%%%%%%%%%
\definition
Let $[\eta]\in H^{1,1}(M)$ be a real (1,1)-class on
a hyperk\"ahler manifold $M$. We say that $[\eta]$
is {\bf parabolic} if $q([\eta],[\eta])=0$.
A line bundle $L$ is called {\bf parabolic} if $c_1(L)$
is parabolic.

\hfill

%The terminology is related to the usual intuition which
%comes from the theory of partial differential equations.
%A K\"ahler form is ``elliptic'' because the corresponding
%Laplacian is an elliptic operator. The Laplacians
%associated with the parabolic forms have symbol which
%is positive semidefinite, that is, parabolic.

\hfill

The present paper is a study of
algebro-geometric properties of parabolic
bundles and cohomology classes, in hope to 
find criteria for effectivity.

%%%%%%%%%%%%%%%%%%%%%%%%%%%%%%%%%%%%%%%%%%%%%%%%%%%%%%%%%%%%
\subsection{The hyperk\"ahler SYZ conjecture}
%%%%%%%%%%%%%%%%%%%%%%%%%%%%%%%%%%%%%%%%%%%%%%%%%%%%%%%%%%%%

\theorem
(D. Matsushita, see \cite{_Matsushita:fibred_}).
Let $\pi:\; M \arrow X$ be a surjective holomorphic map
from a hyperk\"ahler manifold $M$ to $X$, with $0<\dim X < \dim M$.
Then $\dim X = 1/2 \dim M$, and the fibers of $\pi$ are 
holomorphic Lagrangian (this means that the symplectic
form vanishes on the fibers).\footnote{Here, as elsewhere,
we silently assume that the hyperk\"ahler manifold $M$ is
simple.}

\hfill

\definition Such a map is called
{\bf a holomorphic Lagrangian fibration}.

\hfill

\remark The base of $\pi$ is conjectured to be
rational. J.-M. Hwang (\cite{_Hwang:base_}) 
proved that $X\cong \C P^n$, if it is smooth.
D. Matsushita (\cite{_Matsushita:CP^n_}) 
proved that it has the same rational cohomology
as $\C P^n$.

\hfill

\remark
 The base of $\pi$ has a natural flat connection
on the smooth locus of $\pi$. The combinatorics of this connection
can be used to determine the topology of $M$ 
(\cite{_Kontsevich-Soibelman:torus_},  
\cite{_Gross:SYZ_}),

\hfill

\remark 
Matsushita's theorem is implied by the following formula
of Fujiki. Let $M$ be a hyperk\"ahler manifold, $\dim_\C
M=2n$, and $\eta_1, ..., \eta_{2n}\in H^2(M)$ cohomology
classes. Then 
\begin{equation}\label{_Fujiki_multi_Equation_}
\eta_1\wedge \eta_2 \wedge ... = 
\mu \sum_{\sigma}q(\eta_{\sigma_1} \eta_{\sigma_2})
q(\eta_{\sigma_3} \eta_{\sigma_3})q(\eta_{\sigma_{2n-1}} \eta_{\sigma_{2n}})
\end{equation}
with the sum taken over all permutations, and $\mu$ is a rational
constant depending on the dimension $n$.
An algebraic argument (see e.g. \ref{_product_vanishes_Corollary_})
allows to deduce from this formula that
for any non-zero $\eta \in H^2(M)$,
one would have  $\eta^{n}\neq 0$, and $\eta^{n+1}=0$, 
if $q(\eta,\eta)=0$, and $\eta^{2n}\neq 0$ otherwise.
Applying this to the pullback $\pi^*\omega_X$
of the K\"ahler class from $X$, we immediately
obtain that $\dim_\C X=n$ or $\dim_\C X=2n$. Indeed, 
$\omega_X^{\dim_\C X}\neq 0$ and
$\omega_X^{\dim_\C X+1}= 0$.

\hfill

\definition
 Let $(M,\omega)$ be a Calabi-Yau manifold,
$\Omega$ the holomorphic volume form, and $Z\subset M$ a real 
analytic subvariety, Lagrangian with respect to
$\omega$. If $\Omega\restrict Z$ is proportional to
the Riemannian volume form, $Z$ is called {\bf special
Lagrangian} (SpLag).

\hfill

The special Lagrangian varieties were defined
in  \cite{_Harvey_Lawson:Calibrated_}
by Harvey and Lawson, who proved that they
minimize the Riemannian
volume in their cohomology class. This implies, in
particular, that their moduli are finite-dimensional.
In \cite{_McLean:SpLag_}, McLean studied deformations
of non-singular special Lagrangian
subvarieties and showed that they are unobstructed.

In \cite{_SYZ:MS_is_T_du_}, Strominger-Yau-Zaslow 
tried to explain the mirror symmetry phenomenon
using the special Lagrangian fibrations. They
conjectured that any Calabi-Yau manifold admits
a Lagrangian fibration with special Lagrangian fibers.
Taking its dual fibration, one obtains ``the mirror dual''
Calabi-Yau manifold.

\hfill

\remark
It is easy to see that  a holomorphic Lagrangian
subvariety of a hyperk\"ahler manifold $(M,I)$ is special
Lagrangian on $(M,J)$, where $(I,J,K)$ is a quaternionic
structure associated with the hyperk\"ahler structure
on $M$ (Subsection \ref{_gene_hyper_subva_Subsection_}). 
Therefore, existence of holomorphic
Lagrangian fibrations implies existence
of special Lagrangian fibrations postulated
by Strominger-Yau-Zaslow.

\hfill

\definition A line bundle is called
{\bf semiample} if  $L^N$ is generated
by its holomorphic sections, which have 
no common zeros.

\hfill

%%%%%%%%%%%%%%%%%%%%%%%%%%%%%%%%%%%%%%%%%%%%%%%%%%%%%%%
\remark From the semiampleness 
it obviously follows that $L$ is nef. Indeed,
let $\pi:\; M \arrow {\Bbb P}H^0(L^N)^*$ be the the standard
map. Since sections of $L$ have no common zeros, $\pi$ is 
holomorphic. Then $L\cong \pi^* \calo(1)$, and the
curvature of $L$ is a pullback of 
the K\"ahler form on $\C P^n$. However,
the converse is false: a nef bundle is not 
necessarily semiample 
(see e.g. \cite[Example 1.7]{_Demailly_Peternell_Schneider:nef_}).

\hfill

\remark Let $\pi:\; M \arrow X$ 
be a holomorphic Lagrangian fibration, and $\omega_X$
a K\"ahler class on $X$. Then $\eta:=\pi^*\omega_X$ is 
semipositive, and the corresponding line bundle
is semiample and parabolic. The converse is also
true, by Matsushita's theorem:
if $L$ is semiample and parabolic, $L$ induces a Lagrangian
fibration. This is the only known 
source of non-trivial special Lagrangian fibrations
on Calabi-Yau manifolds.

\hfill

%%%%%%%%%%%%%%%%%%%%%%%%%%%%%%%%%%%%%%%%%%%
\conjecture\label{_SYZ_conj_Conjecture_}
(Hyperk\"ahler SYZ conjecture)
Let $L$ be a parabolic nef line bundle
on a hyperk\"ahler manifold. Then
$L$ is semiample.

\hfill

\remark
This conjecture was stated by many people
(Tyurin, Bogomolov, Hassett-Tschinkel,
Huybrechts, Sawon); please see 
\cite{_Sawon_} for an interesting
and historically important discussion,
and \cite{_Verbitsky:SYZ_} 
for details and reference.

\hfill

\remark
The SYZ conjecture can be seen as
a hyperk\"ahler version of ``abundance conjecture''
(see e.g. \cite{_Demailly_Peternell_Schneider:ps-eff_}, 
2.7.2). 

%%%%%%%%%%%%%%%%%%%%%%%%%%%%%%%%%%%%%%%%%%%%%%%%%%%%%%%%%%%%
\subsection{Lelong numbers and hyperk\"ahler geometry}
%%%%%%%%%%%%%%%%%%%%%%%%%%%%%%%%%%%%%%%%%%%%%%%%%%%%%%%%%%%%

In \cite{_Verbitsky:SYZ_}, it was shown that
any parabolic line bundle $L$ with a smooth metric
of semipositive curvature is $\Q$-effective
(this means that $L^n$ is effective, for some 
integer $n>0$). Further results in this
direction require detailed study of
singularities of positive currents
on hyperk\"ahler manifolds.
The present paper is an attempt
to understand these singularities.

\hfill

Let $[\eta]$ be a nef cohomology class.
Using weak compactness of positive currents,
it is possible to show that $[\eta]$ is represented
by a positive, closed $(1,1)$-current $\eta$ 
(\ref{_weak_compa_curre_Claim_}). Locally, $\eta$ can be
considered as a curvature of a 
singular metric on a line bundle.

 Using a local $dd^c$-lemma, we may assume that $\eta=dd^c \phi$,
 for some function $\phi$, which is plurisubharmonic, because
 $\eta$ is positive. Then $\eta$ is a curvature of a 
 trivial bundle with a singular metric $h \arrow e^{-2\phi}|h|^2$.
 
 A {\bf multiplier ideal sheaf} ${\cal I}(\eta)$ of a current $\eta$ is 
 an ideal of all holomorphic functions $h$ on $M$
 for which $e^{-2\phi}|h|^2$ is locally integrable.
 A. M. Nadel has shown that a multiplier ideal sheaf of
 a positive current is always coherent.
 
The notion of a multiplier ideal has
many applications in algebraic geometry, due to the
Nadel's vanishing theorem. 

\hfill

%%%%%%%%%%%%%%%%%%%%%%%%%%%%%%%%%%%%%%%%%%%%%%%%%%%%
\theorem
(Nadel's Vanishing Theorem; see \cite{_Nadel:vanishing_}, 
\cite{_Demailly:vanishing_}). 
Let $(M,\omega)$ be a K\"ahler
manifold, $\eta$ a closed, positive (1,1)-current, $\eta>\epsilon \omega$,
 and $L$ a holomorphic line bundle with $c_1(L)=[\eta]$. 
Consider a singular metric on $L$ associated with
$\eta$, and let ${\cal I}(L)$ be the sheaf of $L^2$-integrable
sections. Then $H^i({\cal I}(L)\otimes K_M)=0$ for all
$i>0$. \endproof

\hfill

The Lelong number $\nu_x(\Theta)$ 
of a $(p,p)$-current $\Theta$ at $x\in M$, as defined in 
\cite{_Demailly:ecole_}, is 
mass of a measure $\Theta \wedge \mu_x^{n-p}$
carried at $x$, where $\mu_x= dd^c(\log \dist_x^2)$, and $\dist_x^2$ is
a square of a distance from $x$. The current
$\mu_x$ can be approximated by smooth, closed, positive
currents using a regularized maximum function (see 
Subsection \ref{_regularization_Subsection_}), and this
allows one to define the product $\Theta \wedge \mu_x^{n-p}$ 
as a limit of closed, positive currents with bounded mass,
well defined because of a weak compactness principle.

 For a positive number $c>0$, the Lelong set $F_c$
 of a (1,1)-current $\eta$ is a set of all points 
 $x\in M$ with $\nu_x(\eta) \geq c$.
By Siu's theorem (\cite{_Siu:1974_}),
a Lelong set of a positive, closed 
current is complex analytic.

\hfill

The following theorem was proven in \cite{_Verbitsky:SYZ_},
using an advanced version of Nadel's vanishing, due to 
\cite{_Demailly_Peternell_Schneider:ps-eff_}.

\hfill

%%%%%%%%%%%%%%%%%%%%%%%%%%%%%%%%%%%%%%%%%%%%%%%
\theorem\label{_Lelong_Q_eff_Theorem_}
(\cite[Theorem 4.1]{_Verbitsky:SYZ_})
 Let $L$ be a parabolic nef bundle
on a hyperk\"ahler manifold, and $\eta$
a positive closed current, representing $c_1(L)$. Assume that all Lelong
numbers of $\eta$ vanish. Then $L$ is $\Q$-effective.
\endproof

\hfill

In the present paper, we show that 
at least one of the Lelong sets of a parabolic nef current
on a hyperk\"ahler manifold is coisotropic
with respect to its holomorphic symplectic
form (\ref{_nef_Lelong_sets_coiso_Corollary_}),
unless all Lelong numbers vanish.

Comparing \ref{_nef_Lelong_sets_coiso_Corollary_}
and \ref{_Lelong_Q_eff_Theorem_},
we obtain the following. Let $L$ be a  
parabolic nef bundle
on a hyperk\"ahler manifold $M$. Then either $L$
is $\Q$-effective, or $M$ has non-trivial
coisotropic subvarieties. A similar result
was proven in by Campana-Oguiso-Peternell, who have shown that
such a manifold always contains a subvariety
of dimension $\geq 2$ (\cite[Theorem 6.2]{_COP:non-alge_}).

For a generic hyperk\"ahler manifold,
all complex subvarieties are holomorphically
symplectic (\cite{_Verbitsky:Symplectic_I_}, 
\cite{_Verbitsky:Symplectic_II_}). Therefore,
such a manifold does not have any 
coisotropic subvarieties
(\ref{_holomo_sy_subva_Remark_}). 
This implies that all Lelong numbers of a parabolic nef current 
on a generic hyperk\"ahler manifold vanish 
(\ref{_Lelong_sets_generic_Corollary_}).

%%%%%%%%%%%%%%%%%%%%%%%%%%%%%%%%%%%%%%%%%%%%%%%%%%%%%%%%%%%%

\section{Hyperk\"ahler geometry: preliminary results}

%%%%%%%%%%%%%%%%%%%%%%%%%%%%%%%%%%%%%%%%%%%%%%%%%%%%%%%%%%%%

%%%%%%%%%%%%%%%%%%%%%%%%%%%%%%%%%%%%%%%%%%%%%%%%%%%%%%%%%%%%
\subsection{The structure of a K\"ahler cone}
%%%%%%%%%%%%%%%%%%%%%%%%%%%%%%%%%%%%%%%%%%%%%%%%%%%%%%%%%%%%

\definition
A class $\eta\in H^{1,1}(M)$ is called {\bf pseudoeffective}
if it can be represented by a positive current, and {\bf nef}
if it lies in the closure of the K\"ahler cone.

\hfill

The following useful theorem, due to S. Boucksom, is known as
{\bf the divisorial Zariski decomposition theorem}.

\hfill

%%%%%%%%%%%%%%%%%%%%%%%%%%%%%%%%%%%%%%%%%%%%%%%%%
\theorem
(\cite{_Boucksom_})
Let $M$ be a  hyperk\"ahler manifold. Then every
pseudoeffective class can be decomposed as a sum
\[
\eta = \nu + \sum_i a_i E_i,
\]
where $\nu$ is nef, $a_i$ positive numbers, and 
$E_i$ exceptional divisors satisfying $q(E_i,E_i)<0$. 
\endproof

\hfill

%%%%%%%%%%%%%%%%%%%%%%%%%%%%%%%%%%%%%%%%%%%%%%%%%%%
\remark 
Let $M_1, M_2$ be holomorphic symplectic
manifolds, bimeromorphically equivalent. Then 
$H^2(M_1)$ is naturally isomorphic to $H^2(M_2)$, and this
isomorphism is compatible with Bogomolov-Beauville-Fujiki form.
Indeed, the manifolds $M_i$ have trivial canonical bundle,
hence a bimeromorphic equivalence is non-singular
in codimension 1. 

\hfill

%%%%%%%%%%%%%%%%%%%%%%%%%%%%%%%%%%%%%%%%%%%%%%%%%%%
\definition
A {\bf modified nef cone} (also ``birational nef cone'' and
``movable nef cone'') is a closure of a union of all nef cones for
all bimeromorphic models of a holomorphically symplectic manifold  $M$.

\hfill

%%%%%%%%%%%%%%%%%%%%%%%%%%%%%%%%%%%%%%%%%%%
\theorem
(\cite{_Boucksom_}, \cite{_Huybrechts:cone_}).
On a hyperk\"ahler manifold,
the modified nef cone is  dual to the pseudoeffective
cone under the Bogomolov-Beauville-Fujiki pairing.
\endproof

\hfill

%%%%%%%%%%%%%%%%%%%%%%%%%%%%%%%%%%%%%%%%%%%
\corollary\label{_Kahler_cone_Pic=1_Corollary_}
 Let $M$ be a simple hyperk\"ahler manifold
such that all integer $(1,1)$-classes satisfy $q(\nu,\nu)\geq 0$.
Then its K\"ahler cone is one of two components $K_+$ of a set 
$K:=\{ \nu\in H^{1,1}(M,\R)\ \ |\ \ q(\nu,\nu)>0\}$.

\hfill

\proof 
The pseudoeffective cone $K_{ps}$ of $M$ is 
equal to the nef cone $K_n$ by the divisorial Zariski decomposition. 
A square of a K\"ahler form is positive,
hence $K_n=K_{ps}$ is contained in one of components
of $K$, denoted by $K_+$. 
This gives inclusions
\begin{equation}\label{_cone_inclu_Equation_}
K_{ps}=K_n\subset K_{mn}\subset K_+
\end{equation}
Since $K_+$ is self-dual,
dualising \eqref{_cone_inclu_Equation_} gives
\begin{equation}\label{_dual_inclu_Equation_}
K_+ \subset K_{ps} \subset K_{mn}=K_n^*
\end{equation}
However, all elements of $K_{mn}$ satisfy $q(\eta,\eta)\geq 0$,
hence $K_{mn}\subset K_+$. Then \eqref{_dual_inclu_Equation_}
gives
\[
K_+ \subset K_{ps} \subset K_{mn}=K_n^*\subset K_+,
\]
and all these cones are equal.
\endproof

\hfill

\remark From the Hodge index theorem, it follows
immediately that the condition 
\[ 
\forall \eta\in \Pic(M)\ \   q(\eta,\eta)\geq 0
\] implies that $\Pic(M)$ has rank 1.

\hfill

%%%%%%%%%%%%%%%%%%%%%%%%%%%%%%%%%%
\remark\label{_Pic=1_Parabolic_Remark_}
From \ref{_Kahler_cone_Pic=1_Corollary_},
it follows that on a hyperk\"ahler manifold
with $\Pic(M)=\Z$, for any rational 
class $\eta\in H^{1,1}(M)$ with $q(\eta,\eta)\geq 0$,
either $\eta$ or $-\eta$ is nef.

%%%%%%%%%%%%%%%%%%%%%%%%%%%%%%%%%%%%%%%%%%%%%%%%%%%%%%%%%%%%
\subsection{Subvarieties in generic hyperk\"ahler
  manifolds}
\label{_gene_hyper_subva_Subsection_}
%%%%%%%%%%%%%%%%%%%%%%%%%%%%%%%%%%%%%%%%%%%%%%%%%%%%%%%%%%%%

This is a brief introduction to the theory
of subvarieties in generic hyperk\"ahler manifolds.
For more details and missing reference, please see
\cite{_Verbitsky:Symplectic_II_} and 
\cite{_Verbitsky:hypercomple_}. 

Recall now that any K\"ahler manifold 
with trivial canonical class admits a unique
Ricci-flat K\"ahler metric in any given K\"ahler class
(\cite{_Yau:Calabi-Yau_}).   Using Bochner's vanishing,
it is possible to show that any holomorphic form
on a compact Ricci-flat manifold is parallel
with respect to the Levi-Civita connection.

If the manifold $M$ is holomorphically symplectic, a Ricci-flat
metric together with the holomorphic symplectic form 
can be used to construct a triple of complex structures
$(I, J, K)$ satisfying quaternionic relations $I\circ J =
- J\circ I =K$, and parallel with respect to the
Levi-Civita connection. In differential geometry and physics,
hyperk\"ahler manifolds are usually defined
in terms of this quaternionic structure
(\cite{_Besse:Einst_Manifo_}).

Consider an operator $L= aI + b J + c K$,
with $a, b, c \in \R$ satisfying $a^2 + b ^2+ c^2 =1$.
Since $I, J, K$ are parallel with respect to 
the Levi-Civita connection, $L$ is also parallel.
Using the quaternionic relations, we obtain $L^2=-1$.
Since $L$ is parallel, it is an integrable complex
structure. Such a complex structure is called 
{\bf induced by the quaternionic action}. 
The set of induced complex structures is parametrized by 
the 2-dimensional sphere $S^2$. It is easy to check that
this gives a holomorphic family of complex
structures on $M$ over $\C P^1$. The total space
of this family is called {\bf the twistor space} of $M$.
Denote the base of the twistor family by $C$, $C \cong \C P^1$.

The group $SU(2)$ of unitary quaternions acts on $TM$.
We extend this action to the bundle $\Lambda^*M$
of differential forms by multiplicativity. 
This action is parallel, hence it commutes
with the Laplacian. This gives a natural $SU(2)$-action
on $H^*(M)$, analogous to the Hodge decomposition
in K\"ahler geometry.

Given a class $v\in H^{2p}(M)$ which is not $SU(2)$-invariant,
let $S_v\subset C$ be the set of all induced complex
structures $L\in C$ for which $v\in H^{p,p}(M)$.
For an $SU(2)$-class, we set $S_v = \emptyset$. 
Since the Hodge decomposition on $(M, L)$
is induced by the $SU(2)$-action, $S_v$ can be expressed
through the action of $SU(2)$. Then it is easy to check that
$S_v$ is finite, for all $v$.

The union $R:= \bigcup_{v\in H^*(M, \Z)} S_v$ is countable.
Clearly, for any induced complex structure $L\notin R$,
\[
v\in H^{p,p}(M) \cap H^{2p}(M, \Z) \Rightarrow 
\text{$v$ is $SU(2)$-invariant}.
\]
%%%%%%%%%%%%%%%%%%%%%%%%%%%%%%%%%%%%%%%%%%%%%%%%
\definition
An induced complex structure $L$ is called {\bf generic}
if $L\notin R$. 

\hfill

As shown in \cite{_Verbitsky:Symplectic_I_},
a closed complex subvariety $X \subset M$ with 
fundamental class $[X]\in H^{2p}(M)$ 
$SU(2)$-invariant is necessarily holomorphically symplectic
outside of its singularities. 

\hfill

%%%%%%%%%%%%%%%%%%%%%%%%%%%%%%%%%%%%%%%%%%%%%%%%
\theorem \label{_gene_subva_Theorem_}
(\cite{_Verbitsky:Symplectic_I_})
Let $(M,I,J,K)$ be a hyperk\"ahler manifold equipped with a
quaternionic structure, and $L$ a generic induced complex
structure. Then all complex subvarieties $X\subset (M, L)$
are holomorphically symplectic outside of singularities.
\endproof

\hfill

\remark
In \cite{_Verbitsky:hypercomple_}
it was also shown that a normalization of $X$ is smooth
and holomorphically symplectic.

\hfill

%%%%%%%%%%%%%%%%%%%%%%%%%%%%%%%%%%%%%%%%%%%%%%%%
\definition
A hyperk\"ahler manifold $(M,I)$ is {\bf generic} if
$I$ is generic for some quaternionic structure
constructed as above.

\hfill

%%%%%%%%%%%%%%%%%%%%%%%%%%%%%%%%%%%%%%%%%%%%%%%%
\remark\label{_holomo_sy_subva_Remark_}
Let $M$ be a generic hyperk\"ahler manifold. Then all
complex subvarieties of $M$ are holomorphically
symplectic, by \ref{_gene_subva_Theorem_}.
In particular, $M$ has no divisors.

%%%%%%%%%%%%%%%%%%%%%%%%%%%%%%%%%%%%%%%%%%%%%%%%%%%%%%%%%%%%
\subsection{Cohomology of hyperk\"ahler manifolds}
%%%%%%%%%%%%%%%%%%%%%%%%%%%%%%%%%%%%%%%%%%%%%%%%%%%%%%%%%%%%

In the sequel, some basic results about cohomology
of hyperk\"ahler manifolds will be used.
The following theorem was proving in \cite{_Verbitsky:cohomo_},
using representation theory.

\hfill

%%%%%%%%%%%%%%%%%%%%%%%%%%%%%%%%%%%%%%%%%%%%%%%%
\theorem \label{_symme_coho_Theorem_}
(\cite{_Verbitsky:cohomo_}) 
Let $M$ be a simple hyperk\"ahler manifold, 
and $H^*_r(M)$ the part of cohomology generated by $H^2(M)$.
Then $H^*_r(M)$ is isomorphic to the symmetric algebra 
(up to the middle degree). Moreover, the Poincare 
pairing on $H^*_r(M)$ is non-degenerate. 
\endproof

\hfill
 
This brings the following corollary.

\hfill

%%%%%%%%%%%%%%%%%%%%%%%%%%%%%%%%%%%%%%%%%%%%%%%%%%%%
\corollary \label{_product_vanishes_Corollary_}
Let $\eta_1, ... \eta_{n+1}\in H^2(M)$ be cohomology
classes on a simple hyperk\"ahler manifold, $\dim_\C M = 2n$.
Suppose that $q(\eta_i, \eta_j)=0$ for all $i, j$. Then
$\eta_1 \wedge \eta_2 \wedge ... \wedge \eta_{n+1}=0$.

\hfill

{\bf Proof:}
Let $H:=\eta_1 \wedge \eta_2 \wedge ... \wedge \eta_{n+1}$.
 From the Fujiki's formula \eqref{_Fujiki_multi_Equation_} it follows
directly that 
\[ 
   H\wedge
   \rho_1 \wedge ... \wedge \rho_{n-1}=0,
\]
for any cohomology classes $\rho_1, ..., \rho_{n-1}\in H^2(M)$.
Therefore, for any $v \in H^{2n-2}_r(M)$,
$H\wedge v=0$. Since the Poincare form is non-degenerate
on $H^{2n-2}_r(M)$ (\ref{_symme_coho_Theorem_}), 
this implies that $H=0$. \endproof

%%%%%%%%%%%%%%%%%%%%%%%%%%%%%%%%%%%%%%%%%%%%%%%%%%%%%%%%%%%%

\section{Cohomology classes dominated by a nef class}

%%%%%%%%%%%%%%%%%%%%%%%%%%%%%%%%%%%%%%%%%%%%%%%%%%%%%%%%%%%%

%%%%%%%%%%%%%%%%%%%%%%%%%%%%%%%%%%%%%%%%%%%%%%%%%%%%%%%%%%%%
\subsection{Positive forms and positive currents}
%%%%%%%%%%%%%%%%%%%%%%%%%%%%%%%%%%%%%%%%%%%%%%%%%%%%%%%%%%%%

In this Subsection, we recall standard 
notions of positivity for $(p,p)$-forms and
currents. A reader may consult 
\cite{_Demailly:ecole_} for more details.

Recall that a real $(p,p)$-form $\eta$
on a complex manifold is called {\bf weakly positive}
if for any complex subspace $V\subset T M$, 
$\dim_\C V=p$, the restriction $\rho\restrict V$
is a non-negative volume form. Equivalently,
this means that 
\[ 
  (\1)^p\rho(x_1, \bar x_1, x_2, \bar x_2, ..., x_p, \bar
  x_p)\geq 0,
\]
for any vectors $x_1, ... x_p\in T_x^{1,0}M$.
A form is called {\bf strongly positive} if it can 
be locally expressed as a sum
\[
\eta = (-\1)^p\sum_{i_1, ... i_p} 
\alpha_{i_1, ... i_p} \xi_{i_1} \wedge \bar\xi_{i_1}\wedge ... 
\wedge \xi_{i_p} \wedge \bar\xi_{i_p}, \ \  
\]
running over some set of $p$-tuples 
$\xi_{i_1}, \xi_{i_2}, ..., \xi_{i_p}\in \Lambda^{1,0}(M)$,
with $\alpha_{i_1, ..., i_p}$ real and non-negative functions on $M$.

The strongly positive and the weakly positive forms
form closed, convex cones in the space 
$\Lambda^{p,p}(M,\R)$ of real $(p,p)$-forms.
These two cones are dual with respect to the Poincare pairing
\[
\Lambda^{p,p}(M,\R) \times \Lambda^{n-p,n-p}(M,\R)\arrow \Lambda^{n,n}(M,\R)
\]
For (1,1)-forms and $(n-1,n-1)$-forms,
the strong positivity is equivalent
to weak positivity.

\hfill

%%%%%%%%%%%%%%%%%%%%%%%%%%%%%%%%%%%%
\remark\label{_multiindex_Remark_}
A  strongly positive form is a
linear combination of products
\[ 
  \alpha (\1)^p z_{i_1}\wedge \bar z_{i_1}\wedge
   z_{i_2}\wedge \bar z_{i_2} \wedge z_{i_k}\wedge \bar z_{i_k}
\]
where $\alpha$ is a smooth, positive function, and
$z_1, ..., z_n \in \Lambda^{1,0}(M)$ is a basis in
$(0,1)$ forms. In the sequel, we shall abbreviate such
a form as $\alpha (z\wedge\bar z)_I$, where $I=(i_1, ..., i_k)$
is a multiindex.

\hfill

A {\bf current} is a form taking values in distributions.
The space of $(p,q)$-currents on $M$ is denoted by $D^{p,q}(M)$.
A {\bf strongly positive current}\footnote{In the present
paper, we shall often omit ``strongly'', because we are 
only interested in strong positivity.}
 is a linear combination
\[
\sum_I \alpha_I (z\wedge\bar z)_I
\]
where $\alpha_I$ are positive, measurable functions,
and the sum is taken over all multi-indices $I$.
An integration current of a
closed complex subvariety is a strongly
positive current. 

Notice that ``strongly positive'' should not be 
confused with ``strictly positive'' (the latter means that
a class belongs to the inner part of a positive cone). 
For instance, 0 is a strongly positive current.

Positivity of a current
$\nu$ is often expressed as $\nu \geq 0$.
If $\nu_1-\nu_2$ is positive, one often writes
$\nu_1 \geq \nu_2$.

It is easy to define the de Rham differential on 
currents, and check that its cohomology coincides
with the de Rham cohomology of a manifold.

\hfill

{\bf Mass} of a positive $(p,p)$-current $\rho$
on a compact $n$-dimensional K\"ahler manifold $(M, \omega)$ is
a number $\int_M \rho \wedge \omega^{n-p}$. This number is
non-negative, and never vanishes, unless $\rho=0$.

\hfill

%%%%%%%%%%%%%%%%%%%%%%%%%%%%%%%%%%%%
\claim\label{_weak_compa_curre_Claim_}
(``weak compactness of positive currents'')
Let $\{\eta_i\}$ be a sequence of positive
$(p,p)$-currents with bounded mass. Then
$\{\eta_i\}$ has a subsequence converging
to a positive current in weak topology.
\hfil \endproof

\hfill

The de Rham differential is by definition
continuous in the topology of currents, and
the projection from closed currents to the de Rham
cohomology also continuous. Then, weak 
compactness implies the following
useful result.

\hfill

%%%%%%%%%%%%%%%%%%%%%%%%%%%%%%%%%%%%
\corollary\label{_limit_posi_cla_Corollary_}
Let $\eta_i \in H^{p,p}(M)$ be a sequence
of cohomology classes represented by
closed, positive currents, and $\eta$
its limit. Then $\eta$ also can be
represented by a closed, positive current.
\hfil \endproof

\hfill

%%%%%%%%%%%%%%%%%%%%%%%%%%%%%%%%%%%%%%%%%%%%%%%%%%%%%%%%%%%%
\definition
A {\bf nef current} is
a positive, closed current, obtained
as a weak limit of strongly positive, closed forms.

\hfill

%%%%%%%%%%%%%%%%%%%%%%%%%%%%%%%%%%%%%%%%%%%%%%%%%%%%%%%%%%%%
\definition
Let $\eta$, $\eta'$
be nef currents. Choose sequences $\{\eta_i\}$, $\{\eta'_i\}$
of closed, strongly positive forms converging to $\eta$, $\eta'$.
Then $\{\eta_i\wedge\eta'_i\}$ is a bounded sequence of
closed, strongly positive forms. From weak compactness
it follows that $\{\eta_i\wedge\eta'_i\}$ has a limit.
We define {\bf a product} $\eta\wedge \eta'$ of nef currents 
as a form which can be obtained as a limit of 
$\{\eta_i\wedge\eta'_i\}$, for some choice of 
sequences $\{\eta_i\}$, $\{\eta'_i\}$. The limit
$\{\eta_i\wedge\eta'_i\}$ is non-unique (see the
example below). However, it is a closed, positive current, which
represents the product of the corresponding
cohomology classes.

\hfill

%%%%%%%%%%%%%%%%%%%%%%%%%%%%%%%%%%%%%%%%%%%%%%%%%
\example \label{_product_many_Example_}
Let $M=\C P^2$. Given a hyperplane $H$, we choose
a sequence of positive, closed (1,1)-forms
$\eta_i(H)$ converging to the current of integration 
$[H]$ of $H$. Suppose that the absolute value
of $\eta_i(H)$ is bounded everywhere by $C_i$, and the
mass of $[H]- \eta_i(H)$ is bounded by $\epsilon_i$.
Let $\alpha$ be a positive (1,1)-current.
Then the mass of $([H]-\eta_i(H))\wedge \alpha$
is bounded by $\epsilon_i \sup|\alpha|$:
\begin{equation}\label{_prod_bound_Equation_}
 \int_{\C P^2} \bigg |([H]-\eta_i(H))\wedge 
 \alpha\bigg|\leq \epsilon_i \sup|\alpha|
\end{equation}
Let now $H, H'$ be two distinct hyperplanes,
and $\eta_i(H)$, $\eta_i(H')$ the sequences of positive,
closed forms approximating $H, H'$ as above. Then 
\eqref{_prod_bound_Equation_} implies that
\begin{equation}\label{_int_of_diff_prod_Equation_}
\int_{\C P^2}\bigg |([H]-\eta_i(H))\wedge \eta_j(H')\bigg| \leq \epsilon_i C_j.
\end{equation}
Choosing a sequence $i_k, j_k$ in such a way that
$\lim\limits_{k\rightarrow \infty} \epsilon_{i_k} C_{j_k}=0$, 
and applying \eqref{_int_of_diff_prod_Equation_},
we obtain that the sequence $\eta_{i_k}(H)\wedge \eta_{j_k}(H')$
has the same limit as $\lim [H]\wedge \eta_j(H')=[p]$,
where $p=H \cap H'$ is a point where $H$ and $H'$ intersect.
Given a sequence $H_l$ of planes converging to $H$,
with $H_l\cap H=p$, and applying the same argument,
we obtain a sequence $\eta_{i_k}(H)\wedge \eta_{j_k}(H_k)$
converging to $[p]$. However, $\eta_{j_k}(H_k)$, for an appropriate
choice of an approximating sequence, clearly
converges to $H$. This gives a sequence of
closed, positive forms $\eta_i, \eta'_i$
converging to $[H]$, and the product 
$\eta_i\wedge\eta'_i$ converges to 
the current of integration $[p]$, associated
with an arbitrary point $p\in H$.

\hfill

However, there are situations when the product of currents
is well defined.

\hfill

%%%%%%%%%%%%%%%%%%%%%%%%%%%%%%%%%%%%%%%%%%%%%%%%%%%%%%%%%%%%
\claim\label{_cprod_currents_Claim_}
Let $\eta_1, ..., \eta_m$ be positive, closed (1,1)-currents
with isolated singularities. Then the nef current 
$\eta_1\wedge \eta_2\wedge ... \wedge\eta_m$ 
is uniquely defined.

\hfill

{\bf Proof:} 
Let $\eta_i(k)$ be a sequence of smooth, closed,
positive (1,1)-forms converging to $\eta_i$, $k =0,1, 2,3,...$
This statement is local, hence we may assume that
$\eta_i= dd^c\phi_i$, for some plurisubharmonic functions
$\phi_i$ with isolated singularities in a discrete set $Z\subset M$,
and $\eta_i(k)= dd^c\phi_i(k)$, where $\phi_i(k)$ are smooth. 
For any compact subset $K$ not intersecting $Z$, chose 
$\phi_i(k)$ in such a way that the
restrictions of $\phi_i(k)$ to $K$ converge: $\lim_k\phi_i(k)=\phi_i$.
Then the limit function $\lim_k\phi_i(k)=\phi_i$ is uniquely determined,
and the product $\eta_1\wedge \eta_2\wedge ... \wedge\eta_m$ 
is uniquely defined by Chern-Levine-Nirenberg theorem 
(\cite[(2.3)]{_Demailly:Trento_}).
\endproof

%%%%%%%%%%%%%%%%%%%%%%%%%%%%%%%%%%%%%%%%%%%%%%%%%
\subsection{Regularization for nef currents}
\label{_regularization_Subsection_}
%%%%%%%%%%%%%%%%%%%%%%%%%%%%%%%%%%%%%%%%%%%%%%%

In \cite{_Demailly_1982_}, the notion
of a {\em regularized maximum} of two functions was
defined. Choose $\epsilon >0$, and let 
$\max_\epsilon:\; \R^2\arrow \R$ be a smooth, convex
function which is monotonous in both arguments
and satisfies $\max_\epsilon(x, y) = \max(x,y)$
whenever $|x-y|>\epsilon$. Then $\max_\epsilon$ is called
{\bf a regularized maximum}. It is easy to show
(\cite{_Demailly_1982_}) that a regularized maximum
of two strictly plurisubharmonic functions is
again strictly plurisubharmonic. Moreover, 
for any smooth form $A$ and $L^1$-functions
$x,y$ which satisfy $A+dd^cx\geq 0$
and $A+dd^cy\geq 0$, one would have
$A+dd^c\max_\epsilon(x,y)\geq 0$.

\hfill

Recall that {\bf an almost plurisubharmonic function}
is a generalized function $f$ which satisfies $dd^cf +A\geq 0$
for some smooth (1,1)-form $A$. Clearly,
almost plurisubharmonic functions  are locally integrable.

\hfill

 The Demailly's Regularization Theorem 
 (\cite{_Demailly:Reg_}, Theorem 1.1, \cite{_Demailly:ecole_}, 21.3) 
 implies that any positive, closed (1,1)-current $T$ on a K\"ahler
 manifold $(M, \omega)$ can be weakly approximated
 by a sequence $T_k$ of closed, real $(1,1)$-currents 
in the same cohomology class satisfying the following assumptions
 
 \begin{description}
 \item[(i)] $T_k + \delta_k \omega\geq 0 $, where 
 $\{\delta_k\}$ is a sequence of real numbers converging to 0.
 \item[(ii)] $T_k$ are smooth outside of a 
 complex analytic subset $Z_k \subset M$, with
 $Z_1 \subset Z_2 \subset ...$
 \item[(iii)] Let $T_0$ be a smooth form cohomologous to $T$. Then
 $T_k = T_0 +dd^c \psi_k$, where 
 $\psi_k$ is a non-increasing sequence of almost
 plurisubharmonic functions converging to an
almost plurisubharmonic $\psi$,
 which satisfies $dd^c \psi+T_0=T$.
 \item[(iv)]  Locally around $Z_k$, the functions $\psi_k$ have logarithmic
 poles, namely \[ \psi_k = \lambda_k \log \sum|g_{k,l}|^2+\tau_k,\]
 where $g_{k,l}$ are holomorphic functions vanishing on $Z_k$,
and $\tau_k$ is smooth.
 \item[(v)] The Lelong numbers $\nu(T_k, x)$ of 
 $T_k$ are non-decreasing in $k$ for any $x\in M$
 and converge to $\nu(T, x)$.
 \end{description}

%%%%%%%%%%%%%%%%%%%%%%%%%%%%%%%%%%%%%%%%%%%%%%%%%%%%
\claim\label{_regulariza_nef_Claim_}
Let $T=\eta$ be a nef $(1,1)$-current. Then
the corresponding approximation
currents $T_k+\delta_k\omega$ of the Demailly's regularization
procedure can be also chosen nef.

\hfill

\proof
Let $T_0$ be a smooth, closed form cohomologous to $\eta$.
 Then $\eta = T_0 + dd^c\psi$, where $\psi=\lim_\downarrow \psi_k$.
Let $\nu_i$ be a sequence
of smooth functions such that the 
form $T_0+dd^c \nu_i+\epsilon_i \omega$ is positive, closed, and weakly 
converges to $\eta=T_0+dd^c\psi$, for $\epsilon_i$ a sequence
of real numbers converging to 0. Such $\{\nu_i\}$
exists, because $\eta$ is nef. Indeed, there exists a sequence
of smooth, positive forms $\eta_i$ converging
to $\eta$, with the cohomology class $[\eta_i] =[\eta]+ [\alpha_i]$,
where $[\alpha_i]\in H^{1,1}(M)$ converging to 0. Choose
smooth, closed representatives $\alpha_i$ with
$\lim_i (\sup |\alpha_i|) =0$, and set $\epsilon_i =\sup |\alpha_i|$.  
Then $\epsilon_i \omega_i + \alpha_i$ is positive.
Choose now $\nu_i$ in such a way that 
$\eta_i= dd^c \nu_i +\alpha_i+T_0$. Then
$dd^c \nu_i +T_0+\epsilon_i \omega> dd^c \nu_i +\alpha_i+T_0$,
hence positive.

Adding constant terms 
if necessary, we may assume that
$\lim \nu_i=\psi$. Fix $k\in \Z^{>0}$.
The function $\mu_i(k):=\max_\epsilon(\nu_i, \psi_k)$
is smooth, because $\psi_k$ is smooth outside of 
its poles. The limit $\lim\limits_{i\rightarrow \infty}
\mu_i(k)$ is equal to $\max_\epsilon(\psi, \psi_k)=
\psi_k$ (the last equation holds
because $\psi\leq \psi_k$). Therefore, $\mu_i(k)$
converges to $\psi_k$.
On the other hand, $T_0+dd^c \nu_i+\epsilon_i \omega$ is
positive, and $T_0+ dd^c\psi_k+\delta_k\omega$ 
is positive by approximation property. From the
properties of a regularized maximum it follows
that $T_0+(\delta_k+\epsilon_k)\omega+ dd^c\mu_i(k)$ is also
positive. We proved that the current
$T_k+(\delta_k+\epsilon_k)\omega=
T_0 + dd^c\psi_k+(\delta_k+\epsilon_k)\omega$
is nef. \endproof

%%%%%%%%%%%%%%%%%%%%%%%%%%%%%%%%%%%%%%%%%%%%%%%%%%%%%%%%%%%%
\subsection{Cohomology classes dominated by a nef current}
%%%%%%%%%%%%%%%%%%%%%%%%%%%%%%%%%%%%%%%%%%%%%%%%%%%%%%%%%%%%

\newcommand{\Ll}{\Lbrack}
\newcommand{\Rr}{\Rbrack}

%\newcommand{\some}{{\text{\small\rm some}}}
%\newcommand{\all}{{\small\rm all}}

%%%%%%%%%%%%%%%%%%%%%%%%%%%%%%%%%%%%%%%%%%%%%%%%
\definition
Let $\eta$ be a current obtained as a limit
of a sequence positive, closed, smooth forms $\{\eta_i\}$. Denote
by $\Ll \eta^p\Rr$ the set of all limits of $\eta_i^p$,
for all sequences $\{\eta_i\}$ converging to $\eta$. 
Notice that the set $\Ll\eta_i\Rr$ is never empty,
by weak compactness of currents.

\hfill

%%%%%%%%%%%%%%%%%%%%%%%%%%%%%%%%%%%%%%%%%%%%%%%%
\remark
The set $\Ll \eta^p\Rr$ can be quite big, as
seen from \ref{_product_many_Example_}.

\hfill

\definition
Let $c$ be a Lelong number of a positive, closed current,
and $F_c$ the corresponding Lelong set. By Siu's theorem,
$F_c$ is complex analytic.
An irreducible component of $F_c$ is called 
{\bf a Lelong component} of $F_c$.

\hfill

%%%%%%%%%%%%%%%%%%%%%%%%%%%%%%%%%%%%%%%%%%%%%%%%
\definition
 A real $(p,p)$-current $\nu$
is said to be {\bf dominated by $\eta$} if for some
$\nu'\in \Ll \eta^p\Rr$ one has $\nu'\geq\epsilon\nu$, 
for some $\epsilon>0$.

\hfill
 
 For an example of a current dominated by a nef current
 $\eta$, we look at the 
Lelong sets of $\eta$. From Demailly's regularization,
 Siu's decomposition theorem and Demailly's version of
the intersection theory,
(\cite{_Siu:1974_}, \cite{_Demailly:ecole_}, \cite{_Demailly:Reg_}),
 the following result can be easily deduced.
 
 \hfill
 
%%%%%%%%%%%%%%%%%%%%%%%%%%%%%%%%%%%%
\theorem\label{_Lelong_dominated_Theorem_}
 Let $\eta$ be a nef current on $M$, 
 and $Z$ a $p$-dimensional 
Lelong component, which is not contained in other 
Lelong components. Denote by $[Z]$ its
integration current. Then $[Z]$ is dominated by $\eta$.
 
 \hfill
 
\proof 
To prove  \ref{_Lelong_dominated_Theorem_},
one needs to produce
a sequence of smooth forms $\tilde \eta_i$ converging to $\eta$, 
in such a way that $Z$ is a Lelong component
of  $\lim_i\tilde \eta_i^p$.

This is done as follows. We approximate $\eta$ by currents
$\eta_i$ with logarithmic singularities at $Z$. Then we prove that
for any smooth approximation $\eta_i(j)$ of $\eta_i$, 
the limit $\rho:=\lim_j \eta_i(j)^p$ would have $Z$ as
one of its Lelong components. By semi-continuity
of Lelong numbers, then, $Z$ is a Lelong component
of  $\lim_i\lim_j\eta_i(j)^p$.

By Siu's theorem (\cite{_Siu:1974_}, 
\cite[2.10]{_Demailly:ecole_}), it follows immediately
that the Lelong sets of any current $\rho\in \Ll\eta^p\Rr$ 
have dimension $\geq p$.
By Siu's decomposition formula (\cite[2.18]{_Demailly:ecole_}),
each current $\rho \in \Ll \eta^p\Rr$ can be written as
\[
\rho = \sum_i c_i [Z_i] + R,
\]
where $R$ is a positive, closed current,
$Z_i$ are all $p$-dimensional components
of the Lelong set of $\rho$, and $c_i=\nu_x(\rho)$
for a generic point $x\in Z_i$. Therefore,
to prove \ref{_Lelong_dominated_Theorem_},
it suffices to show that for any  
irreducible $p$-dimensional component $Z$ of the 
Lelong set of $\eta$, there exists
$\rho\in \Ll \eta^p\Rr$, such that
\begin{equation}\label{_Lelong_mult_Equation_}
\nu_x(\rho) \geq \lambda \nu_x(\eta)^p
\end{equation}
for a generic point $x\in Z$, where $\lambda$ a positive
constant continuously depending on the Lelong numbers of $\eta$.

Using the Demailly's regularization theorem,
\ref{_regulariza_nef_Claim_} and semicontinuity
of Lelong numbers, we find that it suffices to prove
inequality \eqref{_Lelong_mult_Equation_} for
the nef currents with logarithmic singularities approximating $\eta$.
Indeed, $\eta$ is a limit $\lim \eta_i$ of nef currents
with logarithmic singularities (\ref{_regulariza_nef_Claim_}).
For each of these currents, the inequality
\eqref{_Lelong_mult_Equation_} would give
\begin{equation}\label{_inequa_lelo_eta_i_Equation_}
\nu_x(\rho_i) \geq \lambda_i \nu_x(\eta_i)^p
\end{equation}
where $\rho_i\in \Ll \eta^p_i\Rr$.
From Demailly's regularization, the Lelong numbers
of $\eta_i$ converge to the Lelong numbers of $\eta$, hence
$\lim_i \lambda_i\nu_x(\eta_i)^p= \lambda\nu_x(\eta)^p$.
The semicontinuity of Lelong numbers
implies $\lim_i\nu_x(\rho_i) \leq \nu_x(\rho)$.
Therefore, \eqref{_inequa_lelo_eta_i_Equation_} for $\rho_i$ brings
\[
\nu_x(\rho)\geq \lim_i\nu_x(\rho_i)\geq \lim_i 
\lambda_i\nu_x(\eta_i)^p = \lambda\nu_x(\eta)^p,
\]
proving \eqref{_Lelong_mult_Equation_} for $\rho$.

It remains to prove 
\eqref{_Lelong_mult_Equation_} when the singularities of
$\eta$ are logarithmic. We are going to show that
\eqref{_Lelong_mult_Equation_} holds for any 
$\rho\in \Ll \eta^p\Rr$, and $z\in Z$ a generic point.
This statement can be proven in a local situation,
for any current $\eta$ on an open ball.

Locally, we can
always assume that $\eta= dd^c \psi$, for some
plurisubharmonic function $\psi$ with logarithmic
singularities. Clearly, for a general point $z\in Z$,
there exists a submanifold $M'\subset M$
transversally intersecting $Z$ in $z$.
Then $\psi\restrict {M'}$ is a plurisubharmonic
function with an isolated logarithmic singularity at
$z$. In this case, the product $(\eta\restrict{M'})^p$ is uniquely
defined (\ref{_cprod_currents_Claim_}), and we obtain 
$\nu_z(\rho)=\nu_z((\eta\restrict{M'})^p)=\nu_z(\eta)^p$
(see \cite{_Demailly:Trento_}, the chapter on ``generalized Lelong
numbers'', for an equivalent definition of Lelong numbers
for which this statement becomes a tautology).
We proved \eqref{_Lelong_mult_Equation_};
\ref{_Lelong_dominated_Theorem_} follows. \endproof

\subsection{$\eta$-coisotropic subvarieties and cohomology classes}
%%%%%%%%%%%%%%%%%%%%%%%%%%%%%%%%%%%%%%%%%%%%%%%%%%%%%%%%%%%%

%%%%%%%%%%%%%%%%%%%%%%%%%%%%%%%%%%%%%%%%%%%%%%%%
\definition
Let $M$ be a hyperk\"ahler manifold,
 $[\eta]\in H^{1,1}(M)$ a 
para\-bo\-lic nef class on $M$,
and $\eta$ a nef current representing
$[\eta]$. We say that a subvariety
$Z\subset M$ is {\bf $[\eta]$-coisotropic}
if $\eta$ dominates the current 
of integration $[Z]$.

\hfill

\definition
Let $(M, \Omega)$ be a holomorphically symplectic manifold,
$\dim_\C Z=2n$, and $Z\subset M$ a complex subvariety of codimension $p\leq n$.
Then $Z$ is called {\bf coisotropic} if the restriction
$\Omega^{n-p+1}\restrict Z$ vanishes on all smooth points of $Z$.

\hfill

%%%%%%%%%%%%%%%%%%%%%%%%%%%%%%%%%%%%%%%%%%%%%%%%%%%%%%%%%%%%
\remark \label{_coiso_restri_Remark_}
This is equivalent to $\Omega$ having rank $\leq n-p$ on $TZ$
in the smooth points of $Z$, which is the minimal possible
rank for a $2n-p$-dimensional subspace in a $2n$-dimensional
symplectic space.

\hfill

%%%%%%%%%%%%%%%%%%%%%%%%%%%%%%%%%%%%%%%%%%%%%%%%
\proposition\label{_eta_coiso_Proposition_}
Let $M$ be a hyperk\"ahler manifold,
$[\eta]\in H^{1,1}(M)$ a parabolic nef class on $M$,
and $Z\subset M$ an  $[\eta]$-coisotropic subvariety 
of complex codimension $p$. Then 
\begin{description}
\item[(i)]  $p\leq n$,
\item[(ii)] $Z$ is coisotropic with respect to a holomorphic
symplectic form on $M$,  and
\item[(iii)] $[\eta]^{n-p+1}\restrict Z=0$.
\end{description}
\proof 
Since $[\eta]$ is nef, we may chose a representative 
nef current $\eta$, which is a limit of 
positive, closed forms $\{\eta_i\}$. Choose
this sequence in such a way that $\eta_i^k$
converges for all $k>0$, and denote the
respective limits by $\eta^k$.

The current $\eta^{n+1}$ is by definition
positive, and cohomologous to 0, because 
$[\eta]^{n+1}=0$ (\ref{_product_vanishes_Corollary_}).  
The domination of $Z$ by $\eta$ 
means that $\eta^p-c[Z]$ is strongly positive, for some $c>0$. 
Since $\eta^{n+1}=0$, $\eta^p-c[Z]\geq 0$
implies that
\begin{equation}\label{_eta_to_Z_Equation_}
0 = \eta^{n+1} = \eta^p\wedge \eta^{n-p+1} \geq [Z]\wedge
\eta^{n-p+1}
\end{equation}
Choosing a subsequence in $\eta_i$ if necessary,
we may assume that the restriction
$\eta_i^{n-p+1}\restrict Z$
converges to a positive current. Then
\eqref{_eta_to_Z_Equation_} gives that
$\eta_i^{n-p+1}\restrict Z= [Z]\wedge
\eta^{n-p+1}$ vanishes everywhere.
This proves \ref{_eta_coiso_Proposition_} (i) and (iii).

Let $\Omega$ be a holomorphic symplectic
form on $M$. It is easy to check that 
$\Omega^i\wedge \bar\Omega^i$ is weakly positive.
A product of a strongly positive current
and a weakly positive form is weakly positive,
hence the product 
$\eta^p\wedge \Omega^{n-p+1}\wedge \bar\Omega^{n-p+1}$ is
positive. However, this product is cohomologous
to 0, as follows from \ref{_product_vanishes_Corollary_}, 
and therefore
\[ 
\eta^p\wedge \Omega^{n-p+1}\wedge \bar\Omega^{n-p+1}=0
\]
Using the same argument as above, we obtain

\begin{equation}
0 = \eta^p\wedge \Omega^{n-p+1}\wedge \bar\Omega^{n-p+1}  
\geq [Z]\wedge  \Omega^{n-p+1}\wedge \bar\Omega^{n-p+1},
\end{equation}
hence $\Omega^{n-p+1}\wedge \bar\Omega^{n-p+1}$
vanishes on $Z$. Using \ref{_coiso_restri_Remark_},
we obtain that this is equivalent to $Z$ being
coisotropic. We proved \ref{_eta_coiso_Proposition_} (ii).
\endproof

\hfill

As follows from \ref{_holomo_sy_subva_Remark_}, 
on a generic hyperk\"ahler manifold $M$,
all complex subvarieties are holomorphically symplectic.
Then $M$ does not have non-trivial coisotropic subvarieties.
This gives

\hfill

%%%%%%%%%%%%%%%%%%%%%%%%%%%%%%%%%%%%%%%%%%%%%%%%
\corollary \label{_Lelong_sets_generic_Corollary_}
Let $M$ be a generic hyperk\"ahler manifold,
and $[\eta]\in H^{1,1}(M)$ a parabolic nef
class, represented by a positive current
$\eta$. Then all Lelong numbers of $\eta$
vanish.
\endproof

\hfill

Comparing \ref{_eta_coiso_Proposition_},
\ref{_Lelong_dominated_Theorem_}, and \ref{_Lelong_Q_eff_Theorem_}, 
we obtain the 
following.

\hfill

%%%%%%%%%%%%%%%%%%%%%%%%%%%%%%%%%%%%%%%%%%%%%%%%
\corollary \label{_nef_Lelong_sets_coiso_Corollary_}
Let $L$ be a parabolic line bundle
on a hyperk\"ahler manifold, equipped with
a singular metric with positive curvature
current $\eta$, which is nef,
and $Z$ a component of its Lelong set.
Then $Z$ is $\eta$-coisotropic.
In particular, $\dim Z \geq  \frac 1 2 \dim M$,
and $Z$ is coisotropic with respect to the
standard holomorphic symplectic structure on $M$.
Moreover, either $c_1(L)$ is represented by a rational
divisor, or the Lelong sets of $L$ are non-empty.
\endproof

\hfill

Comparing this with \ref{_Pic=1_Parabolic_Remark_}, we obtain

\hfill

%%%%%%%%%%%%%%%%%%%%%%%%%%%%%%%%%%%%%%%%%%%%%%%%
\corollary \label{_coiso_gene_Corollary_}
Let $M$ be a hyperk\"ahler manifold with $\Pic(M)=\Z$,
and $L$ a line bundle generating $\Pic(M)$.
Assume that $q(L,L)=0$. Then $c_1(M)$ can be represented
by a divisor, or $M$ has non-trivial coisotropic subvarieties.

\endproof

\hfill

\remark
Since all divisors are coisotropic, the first 
alternative in the Corollary above 
in fact implies the second one.

\hfill

{\bf Acknowledgements:}
I am grateful to S. Boucksom, J.-P. Demailly,
 D. Kaledin, A. Kuznetsov and M. Paun 
for many valuable discussions. Many thanks to 
Tony Pantev for a useful e-mail exchange.
An early version of this paper was used as
a source of a mini-series of lectures at
a conference ``Holomorphically symplectic
varieties and moduli spaces'', in Lille,
June 2-6, 2009. I am grateful to the organizers
for this opportunity and to the participants
for their insight and many useful comments.

{\small

\hfill

\noindent {\sc Misha Verbitsky\\
Laboratory of Algebraic Geometry, \\
Faculty of Mathematics, NRU HSE,\\
7 Vavilova Str. Moscow, Russia
 }
}

\end{document}